\documentclass{amsart}

\title{The space of stability conditions for $A_3$-quiver}
\author{Takahisa Shiina}
\address{Academic Support Center, Kogakuin University,
Tokyo, Japan}
\email{kt13423@ns.kogakuin.ac.jp}

\usepackage{latexsym,amssymb}
\usepackage{array}
\usepackage[dvipdfmx]{graphicx}
\usepackage{epic,eepic}
\usepackage[all]{xy}
\usepackage[mathscr]{eucal}
\usepackage{mathrsfs}
\def\C{\mathbb{C}}
\def\R{\mathbb{R}}
\def\Z{\mathbb{Z}}

\def\P{\mathcal{P}}

\def\D{\mathcal{D}}
\def\E{\mathcal{E}}

\def\z{\mathcal{Z}}
\def\stab{\mathrm{Stab}}

\def\T{\mathscr{T}}

\def\hom{\mathop{\rm Hom}\nolimits}

\def\mod{\mathop{\rm mod}\nolimits}

\usepackage{amsthm}
\theoremstyle{definition}
\newtheorem{df}{Definition}[section]
\theoremstyle{plain}
\newtheorem{thm}[df]{Theorem}

\newcommand{\typei}[5]{%
\begin{minipage}{28mm}
\hfil
\begin{xy}
 (0,0) *{#1}="A", (10,0) *{#2}="B", (20,0) *{#3}="C",
 \ar "A";"B"_-{#4} \ar "B";"C"_-{#5}
\end{xy}
\end{minipage}
}
\newcommand{\typeii}[6]{%
\begin{minipage}{28mm}
\hfil
\begin{xy}
 (0,0) *{#1}="A", (10,0) *{#2}="B", (20,0) *{#3}="C",
 \ar "A";"B"_-{#4} \ar @/^/ "A";"C"^-{#5} \ar "B";"C"_-{#6}
\end{xy}
\end{minipage}
}
\newcommand{\typeiii}[5]{%
\begin{minipage}{28mm}
\hfil
\begin{xy}
 (0,0) *{#1}="A", (10,0) *{#2}="B", (20,0) *{#3}="C",
 \ar "A";"B"_-{#4} \ar @/^/ "A";"C"^-{#5}
\end{xy}
\end{minipage}
}
\newcommand{\typeiv}[5]{%
\begin{minipage}{28mm}
\hfil
\begin{xy}
 (0,0) *{#1}="A", (10,0) *{#2}="B", (20,0) *{#3}="C",
 \ar @/^/ "A";"C"^-{#4} \ar "B";"C"_-{#5}
\end{xy}
\end{minipage}
}

\begin{document}

\begin{abstract}
 The author studied in \cite{shi} the covering map property of
 the local homeomorphism associating to the space of stability
 conditions over the $n$-Kronecker quiver.
 In this paper, we discuss the covering map property for
 stability conditions over the Dynkin quiver of type $A_3$.
 The local homeomorphism from a connected component of
 stability conditions over $A_3$ to 3-dimensional complex
 vector space becomes a covering map when we restrict it to the
 complement of six codimension one subspaces.
\end{abstract}

\maketitle

\section{Introduction}

T.~Bridgeland introduced the notion of stability conditions on
triangulated categories (\cite{bri}).
The idea comes from Douglas's work on $\pi$-stability for D-branes in
string theory (\cite{douc}).
Bridgeland defined a topology (generalized metric) on the set of all
(locally finite) stability conditions, $\stab(\T)$, on a
triangulated category $\T$.
Each connected component $\Sigma\subset\stab(\T)$ is equipped with a
local homeomorphism $\z$ to a certain topological vector space
$V(\Sigma)$ (\cite[Theorem~1.2]{bri}).
Moreover Bridgeland showed that $\Sigma$ is a
(possibly infinite-dimensional) complex manifold.

It is an important problem to show simply connectivity of $\Sigma$,
connectivity of $\stab(\T)$, or (universal) covering map property
of the local homeomorphism $\z:\Sigma\rightarrow V(\Sigma)$.
It have been studied intensively, for example,
stability conditions over K3 surfaces, over curves,
and over Kleinian singularities are studied
\cite{brib,bric,bt,iuu,mac,st,tho}.

In \cite{shi}, the author studied the space of stability conditions
$\stab(P_n)$ constructed on the bounded derived category $D^b(P_n)$
of the finite dimensional representations over $n$-Kronecker quiver
$P_n$ -- the quiver with two vertices and
$n$-parallel arrows.
More precisely, we discussed the covering map property of
the local homeomorphism
\[
 \z:\stab(P_n)\longrightarrow\hom_\Z(K(P_n),\C).
\]
As we only refer to the $1$-Kronecker quiver $P_1$,
i.e. the Dynkin quiver of type $A_2$, the restriction of $\z$
onto the complement of three lines,
\[
 \z|_{\z^{-1}(X)}:\z^{-1}(X) \longrightarrow
 \hom_\Z(K(P_1),\C)\setminus(L_0\cup L_1\cup L_2)=X,
\]
is a covering map. Here
$L_i=\{\,Z\in\hom_\Z(K(P_1),\C)\,|\,Z(E_i)=0\,\}$ and
$E_0$, $E_1$, $E_2$ are well-known exceptional objects
of $D^b(P_1)$.

In this paper, we discuss the covering map property
for stability conditions, $\stab(A_3)$,
over the Dynkin quiver of type $A_3$.
Following is the Main Theorem:
\begin{thm}
 \label{thm:main}
 Let $\Sigma(A_3)$ be a connected component of $\stab(A_3)$
 associating to exceptional collections of $D^b(A_3)$.
 The restriction of the local homeomorphism
 $\z:\Sigma(A_3)\rightarrow\hom_\Z(K(A_3),\C)$
 becomes a covering map when we remove, at least, six
 codimension one subspaces from base space, i.e.
 \[
  \z:\z^{-1}\longrightarrow
  \hom_\Z(K(A_3),\C)\setminus
  (L_1\cup L_2\cup L_3\cup L_4\cup L_5\cup L_6)
 \]
 is a covering map,
 where $L_i=\{\,Z\in\hom_\Z(K(A_3),\C)\,|\,Z(E_i)=0\,\}$
 and $E_i$'s are well-known exceptional objects of $D^b(A_3)$.
\end{thm}

The organization of this paper is as follows:
In Section~2 we prepare basic definitions of stability conditions
by Bridgeland (\ref{re:bri}) and useful results by Macr{\`{\i}}
relating to stability conditions on triangulated categories
generate by finitely many exceptional objects (\ref{re:mac}).
Then we prove the Main Theorem in Section~3, including basic
properties of $D^b(A_3)$ and connectivity of $\Sigma(A_3)$.

\subsection*{Acknowledgement}
I would like to thank Dai Tamaki and So Okada for
invaluable advice.

\section{Preliminaries}

In this section we recall definitions and properties for
Bridgeland's stability conditions on triangulated categories,
and recall Macr{\`{\i}}'s works for the space of stability conditions
on triangulated categories generated by exceptional collection.
See \cite{bri} and \cite{mac} for more details.

\subsection{Stability conditions}\label{re:bri}

Let $k$ be a field and $\T$ be a $k$-linear triangulated category.
The \emph{Grothendieck group} of $\T$, $K(\T)$, is the quotient
group of the free abelian group generated by all isomorphism classes
of objects in $\T$ modulo the subgroup generated by the elements of
the form $[A]+[B]-[C]$ for each distinguished triangle
$A \rightarrow C \rightarrow B$ in $\T$.

A \emph{stability condition} on $\T$ consists of
a \emph{central charge} and a \emph{slicing}.
A \emph{central charge} is a group homomorphism $Z:K(\T) \rightarrow \C$.
A \emph{slicing} $\P$ is a family of full additive subcategories $\P(\phi)$
of $\T$ indexed by real numbers $\phi$ satisfying that
$\hom_\T(A_1,A_2)=0$ if $A_i\in\P(\phi_i)$ and $\phi_1>\phi_2$,
$\P(\phi+1)=\P(\phi)[1]$ for all $\phi\in\R$ and
for each nonzero object $X$ there are sequence of maps
$X_0\rightarrow X_1\rightarrow\dots\rightarrow X_n$ in $\T$ and
sequence of real numbers $\phi_1>\dots>\phi_n$ such that
$X_0=0$, $X_n=X$ and $A_i\in\P(\phi_i)$ for all $i$ which are
objects fitting into the distinguished triangle
$X_{i-1}\rightarrow X_i\rightarrow A_i$.
A pair of a central charge and a slicing, $\sigma=(Z,\P)$, called
a \emph{stability condition} if $Z(A)=m(A)\exp(i\pi\phi)$
for any nonzero object $A\in\P(\phi)$ and some $m(A)>0$.
A nonzero object of $\P(\phi)$ is called \emph{semistable of phase}
$\phi$ and a simple object of $\P(\phi)$ \emph{stable}.
A stability condition $(Z,\P)$ is called \emph{locally finite} if there
exists $\varepsilon>0$ such that $\P(\phi-\varepsilon,\phi+\varepsilon)$
is Artinian and Noetherian.

Let $\stab(\T)$ be the set of all locally-finite stability conditions on
$\T$.
A generalized metric on $\stab(\T)$ is defined by
\[
 d(\sigma_1,\sigma_2)=\sup_{0\neq E\in\T}\left\{
 |\phi_{\sigma_2}^-(E)-\phi_{\sigma_1}^-(E)|,
 |\phi_{\sigma_2}^+(E)-\phi_{\sigma_1}^+(E)|,
 \left|\log \frac{m_{\sigma_2}(E)}{m_{\sigma_1}(E)}\right|
 \right\}\in[0,\infty]
\]
for $\sigma_1,\sigma_2\in\stab(\T)$.
Here $\phi_\sigma^-(E)$ is the lowest number $\phi_n$ and
$\phi_\sigma^+(E)$ is the greatest number $\phi_1$ in the sequence
of $E$ associated to $\sigma$ and
$m_\sigma(E)=\sum_{i=1}^n\left|Z(A_i)\right|$.
When we equip a well-defined linear topology on $\hom_\Z(K(\T),\C)$,
Bridgeland showed that there is a natural local homeomorphism.
\begin{thm}
 \cite[Theorem 1.2]{bri}
 For each connected component $\Sigma\subset\stab(\T)$ there are a
 linear subspace $V(\Sigma)\subset\hom_\Z(K(\T),\C)$ and a local
 homeomorphism $\z:\Sigma\rightarrow V(\Sigma)$ which maps
 a stability condition $(Z,\P)$ to its central charge $Z$.
\end{thm}

The additional important structure on $\stab(\T)$ is the right action
of $\tilde{\mathrm{GL}}^+(2,\R)$, the universal covering of
$\mathrm{GL}^+(2,\R)$, and the left action of the autoequivalences
of $\T$ (see \cite{bri}).

\subsection{Triangulated categories generated by exceptional objects}
\label{re:mac}

An object $E$ in $\T$ is \emph{exceptional} if
$\hom^k_\T(E,E)=\C$ if $k=0$ and $=0$ otherwise.
For two exceptional objects $E$ and $F$,
we write $\mathcal{L}_EF$ to be a \emph{left mutation of $F$ by $E$}
and $\mathcal{R}_FE$ a \emph{right mutation of $E$ by $F$},
which are objects fitting into distinguished triangles;
$\mathcal{L}_EF \rightarrow \hom^\bullet(E,F) \otimes E%
\rightarrow F$ and
$E \rightarrow \hom^\bullet(E,F)^* \otimes F%
\rightarrow \mathcal{R}_FE$.

An \emph{exceptional collection} is a sequence
$\E=(E_1,E_2,\dots,E_n)$ of exceptional objects
such that $\hom^k_\T(E_i,E_j)=0$ for all $k$ and $i>j$.
An exceptional collection is called \emph{complete} if $\{E_i\}$
generates $\T$ by shifts and extensions and is called
\emph{Ext-exceptional} if $\hom^{\leq 0}(E_i,E_j)=0$
for all $i\neq j$.
A \emph{left mutation} $\mathcal{L}_i\mathcal{E}$ and
a \emph{right mutation} $\mathcal{R}_i\mathcal{E}$
($1\leq i\leq n-1$) of $\E$ are defined by
\begin{align*}
 \mathcal{L}_i\E
 &=(E_1,\dots,E_{i-1}\mathcal{L}_{E_i}E_{i+1},E_i,E_{i+2},\dots,E_n)
 \\
 \text{and }\,\mathcal{R}_i\E
 &=(E_1,\dots,E_{i-1},E_{i+1},\mathcal{R}_{E_{i+1}}E_{i},E_{i+2},\dots,E_n).
\end{align*}
The mutation of a (complete) exceptional collection becomes again
a (complete) exceptional collection.
The operations $\mathcal{L}_i$ and $\mathcal{R}_i$ are invertible
each other;
$\mathcal{L}_i\mathcal{R}_i=\mathcal{R}_i\mathcal{L}_i=\text{id}$
for each $i$, and they satisfy the braid relation;
$\mathcal{R}_i\mathcal{R}_{i+1}\mathcal{R}_i=%
\mathcal{R}_{i+1}\mathcal{R}_i\mathcal{R}_{i+1}$
and $\mathcal{L}_i\mathcal{L}_{i+1}\mathcal{L}_i%
=\mathcal{L}_{i+1}\mathcal{L}_i\mathcal{L}_{i+1}$.

Macr{\`{\i}} showed how to construct a stability condition
from complete exceptional collections.
The key rules are following two theorem;
To give a stability condition is equivalent to giving a bounded t-structure
and a stability function on its heart \cite[Proposition 5.3]{bri};
Taking a sequence of integers $p=(p_1,\dots,p_n)$ with
$\E_p=(E_1[p_1],\dots,E_n[p_n])$ is an Ext-exceptional collection,
the smallest extension-closed full subcategory generated by $\E_p$, $Q_p$,
is a heart of a bounded t-structure \cite[Lemma 3.2]{mac}.
By fixing $n$-points $z_1,\dots,z_n$ in
$\{\,m\exp(i\pi\phi)\,|\,m>0,\,0<\phi\leq 1\,\}$
and defining $Z_p:K(\T)\rightarrow\C$ by $Z_p(E_i[p_i])=z_i$,
the pair $(Z_p,Q_p)$ occurs a stability condition.

We write $\Theta_\E$ to be the set of all stability conditions defined
from $\E$ via above process up to the action of $\tilde{\mathrm{GL}}^+(2,\R)$.
Macr{\`{\i}} proved that $\Theta_\E$ is homeomorphic to
\[
 \left\{\,(m_1,\dots,m_n,\phi_1,\dots,\phi_n)\in\R^{2n}\,\left|\,%
 \text{$m_i > 0$ for all $i$ and $\phi_i<\phi_j+\alpha_{i,j}$ for $i<j$}%
 \,\right.\right\}
\]
where
\[
 \alpha_{i,j}=\min_{i<\ell_1<\ell_2<\dots<\ell_s<j}%
 \{\, k_{i,\ell_1}+k_{\ell_1,\ell_2}+\dots+k_{\ell_s,j}-s\,\}
\]
and
\[
 k_{i,j}=
 \begin{cases}
  +\infty & \text{if $\hom^k(E_i,E_j)=0$ for all $k$,} \\
  \min\left\{k\,\left|\, \hom^k(E_i,E_j) \neq 0 \right.\right\}
  & \text{otherwise.}
 \end{cases}
\]
A stability condition $\sigma=(Z,\P)$ in $\Theta_\E$
corresponds to $(m_1,\dots,m_n,\phi_1,\dots,\phi_n)$ by
$m_i=|Z(E_i)|$ and $\phi_i=\phi_\sigma(E_i)$.
Moreover, all $E_i$'s are stable in $\sigma\in\Theta_\E$
(\cite[Lemma~3.16]{mac}) and $\Theta_\E$ is an open, connected and
simply connected $(n+1)$-dimensional submanifold
(\cite[Lemma~3.19]{mac}).

The union of the open subsets $\Theta_\mathcal{F}$ over all
iterated mutations $\mathcal{F}$ of $\E$ is denoted as $\Sigma_\E$.
It is also an open and connected $(n+1)$-dimensional submanifold
(\cite[Corollary~3.20]{mac}).
When all exceptional collections of $\T$ can be obtained,
up to shifts, by iterated mutations of a single exceptional collection,
the previous open subset is denoted as $\Sigma(\T)$ and such $\T$
is called {\em constructible}.

\section{Proof of the Main Theorem}

First we recall some basic facts of the Dynkin quiver of type
$A_3$; $1\rightarrow 2\rightarrow 3$.
We denote by $\C A_3$ the path algebra of $A_3$ over the ground
field $\C$, denote by $\mod\C A_3$ the category of finitely
generated $\C A_3$-modules, denote by $D^b(A_3)$ the bounded derived 
category of $\mod\C A_3$ and denote by $\stab(A_3)$ the space of
stability conditions on $D^b(A_3)$.
Note that the category of finitely generated modules over a path
algebra of a quiver is equivalent to the category of finite
dimensional representations over the one.

The Auslander-Reiten quiver of $\mod\C A_3$ is well-known, and
in this paper we write it as follows;
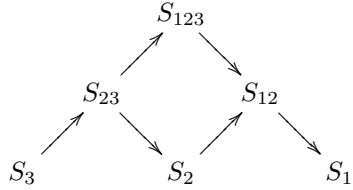
\[
 \xymatrix @!=0.5pc{
   && S_{123} \ar[rd] \\
   & S_{23} \ar[ru] \ar[rd] && S_{12} \ar[rd] \\
   S_3 \ar[ru] && S_2 \ar[ru] && S_1
 }
\]
where $S_1$, $S_2$ and $S_3$ are simple objects relating to
vertices of $A_3$.

The Grothendieck group, $K(A_3)$, of $D^b(A_3)$ is a free abelian
group generated by isomorphism classes of $S_1$, $S_2$ and $S_3$.

The set of all equivalent classes of complete exceptional collections
on $D^b(A_3)$ is shown in \cite{ara} in which two exceptional
collections $(E_1,\dots,E_n)$ and $(F_1,\dots,F_n)$ are called
equivalent if there exists a permutation $\sigma$ and integers
$l_1,\dots,l_n$ such that $F_i=E_{\sigma(i)}[l_i]$ for every $i$.
There are $12$ representatives of complete exceptional collections
on $D^b(A_3)$, listed in Table \ref{listofexc} in our notation.
In the Table, each components
\[
 (X)\, \typeii{S_x}{S_y}{S_z}{a}{b}{c}
\]
means that $X$ is the exceptional collection $(S_x,S_y,S_z)$ such that
\begin{align*}
 a&=\min\left\{\,k\,\left|\,%
   \hom_{D^b(A_3)}^k(S_x,S_y)\neq0\,\right.\right\}, \\
 b&=\min\left\{\,k\,\left|\,%
   \hom_{D^b(A_3)}^k(S_x,S_z)\neq0\,\right.\right\}\,
   \text{ and }\\
 c&=\min\left\{\,k\,\left|\,%
   \hom_{D^b(A_3)}^k(S_y,S_z)\neq0\,\right.\right\},
\end{align*}
however we omit the arrow if $\hom_{\D^b(A_3)}^k(S_x,S_y)=0$
for all $k$.
\begin{table}[h]
\caption{Representatives of exceptional collections on $D^b(A_3)$}
\label{listofexc}
 \begin{tabular}{*{3}{>{$}c<{$}}}
  \hline
  (A)\,\typei{S_1}{S_2}{S_3}{1}{1} &
  (B)\,\typeii{S_2}{S_{12}}{S_3}{0}{1}{1} &
  (I)\,\typeiii{S_{12}}{S_1}{S_3}{0}{1} \\
  (C)\,\typei{S_2}{S_3}{S_{123}}{1}{0} &
  (D)\,\typeii{S_3}{S_{23}}{S_{123}}{0}{0}{0} &
  (J)\,\typeiii{S_{23}}{S_2}{S_{123}}{0}{0} \\
  (E)\,\typei{S_3}{S_{123}}{S_1}{0}{0} &
  (F)\,\typeii{S_{123}}{S_{12}}{S_1}{0}{0}{0} &
  (K)\,\typeiv{S_2}{S_{123}}{S_{12}}{0}{0} \\
  (G)\,\typei{S_{123}}{S_1}{S_2}{0}{1} &
  (H)\,\typeii{S_1}{S_{23}}{S_2}{1}{1}{0} &
  (L)\,\typeiv{S_3}{S_1}{S_{23}}{0}{1} \\\hline
 \end{tabular}
\end{table}

Note that $\hom_{D^b(A_3)}^k(S_x,S_y)\neq 0$ for at most one
$k$ and all $x,y$, and that
$\hom_{D^b(A_3)}^1(S_{12},S_{23})=\C$ and
$\hom_{D^b(A_3)}^0(S_{23},S_{12})=\C$.

\begin{figure}
\[
 \begin{xy}
  (0,50) *+{A}="A",
  (-35,35) *+{B}="B",
  (-50,0) *+{C}="C",
  (-35,-35) *+{D}="D",
  (0,-50) *+{E}="E",
  (35,-35) *+{F}="F",
  (50,0) *+{G}="G",
  (35,35) *+{H}="H",
  (-20,10) *+{I}="I",
  (-20,-10) *+{J}="J",
  (20,10) *+{K}="K",
  (20,-10) *+{L}="L",
  \ar "A";"B"%
  \ar "C";"D"%
  \ar "E";"F"%
  \ar "G";"H"
  \ar "B";"I"%
  \ar "I";"A"%
  \ar "D";"J"%
  \ar "J";"C"%
  \ar "F";"I"%
  \ar "I";"E"%
  \ar "H";"J"%
  \ar "J";"G"%
  \ar @{..>} "B";"C"%
  \ar @{..>} "D";"E"%
  \ar @{..>} "F";"G"%
  \ar @{..>} "H";"A"%
  \ar @{..>} "C";"K"%
  \ar @{..>} "K";"B"%
  \ar @{..>} "E";"L"%
  \ar @{..>} "L";"D"%
  \ar @{..>} "G";"K"%
  \ar @{..>} "K";"F"%
  \ar @{..>} "A";"L"%
  \ar @{..>} "L";"H"%
 \end{xy}
\]
\caption{Mutations of exceptional collections on $D^b(A_3)$}
\label{maingraph}
\end{figure}
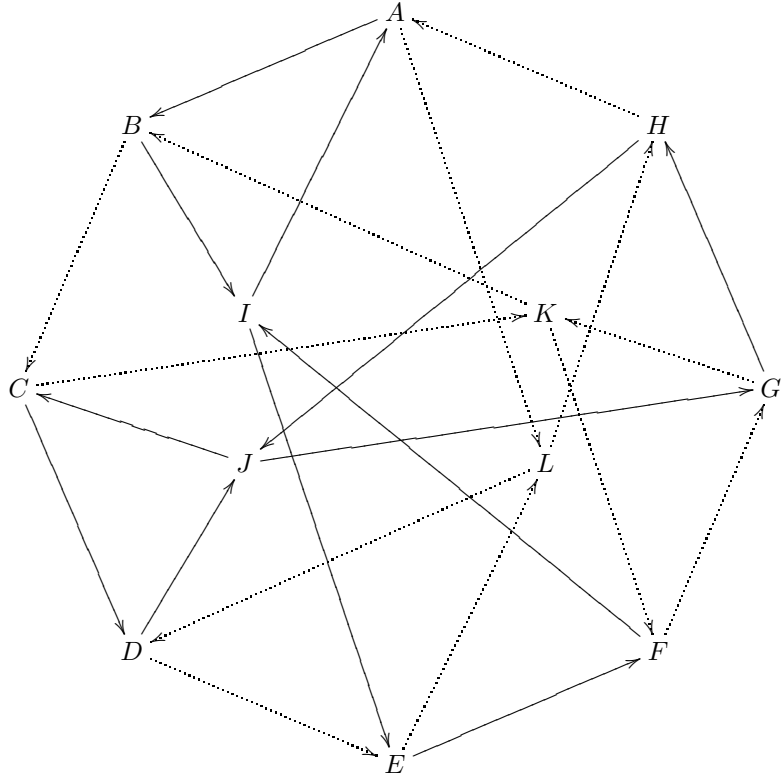

Any two exceptional collections of $D^b(A_3)$ are transitive
by iterated mutations, as illustrated in Figure \ref{maingraph}
in which solid arrows mean $\mathcal{R}_1$, the right mutation
between left and center, and dotted arrows mean $\mathcal{R}_2$,
between center and right.
Applying $\mathcal{R}_2$ (resp. $\mathcal{R}_1$) to $I$ or $J$
(resp. $K$ or $L$) occurs an exceptional collection equivalent
to itself.

\subsection{Applying Macr{\`{\i}}'s method to $D\b(A_3)$}

Let $X=(E_1,E_2,E_3)$ be one of the exceptional collection in
the following list:
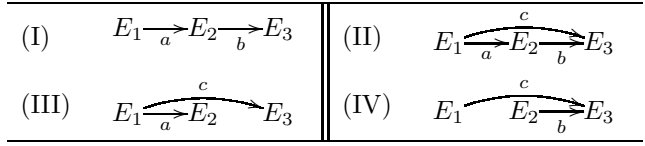
\begin{center}
 \begin{tabular}{l>{$}c<{$}||l>{$}c<{$}}
  \hline
  (I) & \typei{E_1}{E_2}{E_3}{a}{b} &
  (II) & \typeii{E_1}{E_2}{E_3}{a}{c}{b} \\
  (III) & \typeiii{E_1}{E_2}{E_3}{a}{c} &
  (IV) & \typeiv{E_1}{E_2}{E_3}{c}{b} \\\hline
 \end{tabular}
\end{center}
According to Macr{\`{\i}}, $\Theta_\E$ is homeomorphic to
\[
 \left\{\,(m_1,m_2,m_3,\phi_1,\phi_2,\phi_3)\in\R^6\,\left|\,
 \text{$m_i>0$ and $\phi_i<\phi_j+\alpha_{i,j}$ ($i<j$)}
 \,\right.\right\},
\]
where $\alpha_{i,j}$'s and $k_{i,j}$'s are defined as follows:
\begin{center}
 \begin{tabular}{*{5}{>{$}c<{$}}}
  \hline
  & \text{I} & \text{II} & \text{III} & \text{IV} \\\hline
  k_{1,2} & a & a & a & +\infty \\
  k_{2,3} & b & b & +\infty & b \\
  k_{1,3} & +\infty & c & c & c \\\hline
  \alpha_{1,2} & a & a & a & +\infty \\
  \alpha_{2,3} & b & b & +\infty & b \\
  \alpha_{1,3} & a+b-1 & \min\{\, a+b-1,\,c\,\} & c & c
  \\\hline
 \end{tabular}
\end{center}

Applying this to exceptional collections in Table \ref{listofexc},
we obtain relations of $\phi_i$ for each exceptional collection
$X$ as Table \ref{relofce} in which we simply denote $S_x$ instead
of its phase $\phi_i$.
And underlying figures illustrate these relations graphically,
where solid, dashed or dotted arrows mean $Z(S_x)$ in the complex
plane and curved arrows mean that the phase of $S_x$
on arrow-tail is lower than the one of $S_y$ on arrow-head.
\begin{table}[h]
\caption{$\phi_i<\phi_j+\alpha_{i,j}$ for $\Theta_X$}
\label{relofce}
\begin{center}
 \begin{tabular}{>{$}c<{$} *{3}{>{$}l<{$}}}
  \hline
  X & \phi_1,\phi_2 & \phi_2,\phi_3 & \phi_1,\phi_3 \\\hline
  A & S_1<S_2+1 & S_2<S_3+1 & S_1<S_3+1 \\
  B & S_2<S_{12} & S_{12}<S_3+1 & S_2<S_3 \\
  C & S_2<S_3+1 & S_3<S_{123} & S_2<S_{123} \\
  D & S_3<S_{23} & S_{23}<S_{123} & S_3<S_{123}-1 \\
  E & S_3<S_{123} & S_{123}<S_1 & S_3<S_1-1 \\
  F & S_{123}<S_{12} & S_{12}<S_1 & S_{123}<S_1-1 \\
  G & S_{123}<S_1 & S_1<S_{2}+1 & S_{123}<S_2 \\
  H & S_1<S_{23}+1 & S_{23}<S_2 & S_1<S_2 \\
  I & S_{12}<S_1 & & S_{12}<S_3+1 \\
  J & S_{23}<S_2 & & S_{23}<S_{123} \\
  K & & S_{123}<S_{12} & S_2<S_{12} \\
  L & & S_1<S_{23}+1 & S_3<S_{23} \\\hline
 \end{tabular}
\end{center}
\tiny
\begin{tabular}{|c|c|c|c|}
\hline
$A=(S_1,S_2,S_3)$
& $B=(S_2,S_{12},S_3)$
& $C=(S_2,S_3,S_{123})$
& $D=(S_3,S_{23},S_{123})$
\\
\begin{picture}(80,80)(-40,-40)
 \dottedline{3}(-20,-20)(0,0)(20,-20)%
 \put(0,0){\vector(-1,1){20}} \put(-20,20){\makebox(0,0)[rb]{$S_1$}}
 \dashline[50]{5}(0,0)(20,20) \put(20,20){\vector(1,1){0}}
 \put(20,20){\makebox(0,0)[lb]{$S_2$}}
 \dottedline{2}(0,0)(30,10) \put(30,10){\vector(3,1){0}}
 \put(30,10){\makebox(0,0)[l]{$S_3$}}
 \put(20,-20){\makebox(0,0)[lt]{$[-1]$}}
 \put(-20,-20){\makebox(0,0)[rt]{$[-1]$}}
 \arc{28.284}{5.961}{7.069} \put(13.416,4.472){\vector(-1,3){0}}
 \arc{42.426}{5.961}{8.639} \put(20.125,6.708){\vector(-1,3){0}}
 \arc{14.142}{5.498}{7.069} \put(5,5){\vector(-1,1){0}}
\end{picture}
&
\begin{picture}(80,80)(-40,-40)
 \dottedline{3}(-20,-20)(0,0)(20,-20)%
 \put(0,0){\vector(1,1){20}} \put(20,20){\makebox(0,0)[lb]{$S_2$}}
 \dashline[50]{5}(0,0)(-20,20) \put(-20,20){\vector(-1,1){0}}
 \put(-20,20){\makebox(0,0)[rb]{$S_{12}$}}
 \dottedline{2}(0,0)(-10,30) \put(-10,30){\vector(-1,3){0}}
 \put(-10,30){\makebox(0,0)[b]{$S_3$}}
 \put(20,-20){\makebox(0,0)[lt]{$[-1]$}}
 \arc{28.284}{4.390}{5.498} \put(-4.472,13.416){\vector(-3,-1){0}}
 \arc{42.426}{4.390}{7.069} \put(-6.708,20.125){\vector(-3,-1){0}}
 \arc{14.142}{3.927}{5.498} \put(-5,5){\vector(-1,-1){0}}
\end{picture}
&
\begin{picture}(80,80)(-40,-40)
 \dottedline{3}(-20,-20)(0,0)(20,-20)%
 \put(0,0){\vector(-1,1){20}} \put(-20,20){\makebox(0,0)[rb]{$S_2$}}
 \dashline[50]{5}(0,0)(20,20) \put(20,20){\vector(1,1){0}}
 \put(20,20){\makebox(0,0)[lb]{$S_3$}}
 \dottedline{2}(0,0)(-30,-10) \put(-30,-10){\vector(-3,-1){0}}
 \put(-30,-10){\makebox(0,0)[r]{$S_{123}$}}
 \put(20,-20){\makebox(0,0)[lt]{$[-1]$}}
 \arc{28.284}{2.820}{3.927} \put(-13.416,-4.472){\vector(1,-3){0}}
 \arc{42.426}{2.820}{5.498} \put(-20.125,-6.708){\vector(1,-3){0}}
 \arc{14.142}{5.498}{7.069} \put(5,5){\vector(-1,1){0}}
\end{picture}
&
\begin{picture}(80,80)(-40,-40)
 \dottedline{3}(-20,-20)(0,0)(20,-20)%
 \put(0,0){\vector(1,1){20}} \put(20,20){\makebox(0,0)[lb]{$S_3$}}
 \dashline[50]{5}(0,0)(-20,20) \put(-20,20){\vector(-1,1){0}}
 \put(-20,20){\makebox(0,0)[rb]{$S_{23}$}}
 \dottedline{2}(0,0)(10,-30) \put(10,-30){\vector(1,-3){0}}
 \put(10,-30){\makebox(0,0)[t]{$S_{123}$}}
 \put(-20,-20){\makebox(0,0)[rt]{$[+1]$}}
 \arc{28.284}{1.249}{2.356} \put(4.472,-13.416){\vector(3,1){0}}
 \arc{42.426}{1.249}{3.927} \put(6.708,-20.125){\vector(3,1){0}}
 \arc{14.142}{3.927}{5.498} \put(-5,5){\vector(-1,-1){0}}
\end{picture}
\\ \hline
$E=(S_3,S_{123},S_1)$
& $F=(S_{123},S_{12},S_1)$
& $G=(S_{123},S_1,S_2)$
& $H=(S_1,S_{23},S_2)$
\\
\begin{picture}(80,80)(-40,-40)
 \dottedline{3}(-20,-20)(0,0)(20,-20)%
 \put(0,0){\vector(1,1){20}} \put(20,20){\makebox(0,0)[lb]{$S_{3}$}}
 \dashline[50]{5}(0,0)(-20,20) \put(-20,20){\vector(-1,1){0}}
 \put(-20,20){\makebox(0,0)[rb]{$S_{123}$}}
 \dottedline{2}(0,0)(10,-30) \put(10,-30){\vector(1,-3){0}}
 \put(10,-30){\makebox(0,0)[t]{$S_{1}$}}
 \put(-20,-20){\makebox(0,0)[rt]{$[+1]$}}
 \arc{28.284}{1.249}{2.356} \put(4.472,-13.416){\vector(3,1){0}}
 \arc{42.426}{1.249}{3.927} \put(6.708,-20.125){\vector(3,1){0}}
 \arc{14.142}{3.927}{5.498} \put(-5,5){\vector(-1,-1){0}}
\end{picture}
&
\begin{picture}(80,80)(-40,-40)
 \dottedline{3}(-20,-20)(0,0)(20,-20)%
 \put(0,0){\vector(1,1){20}} \put(20,20){\makebox(0,0)[lb]{$S_{123}$}}
 \dashline[50]{5}(0,0)(-20,20) \put(-20,20){\vector(-1,1){0}}
 \put(-20,20){\makebox(0,0)[rb]{$S_{12}$}}
 \dottedline{2}(0,0)(10,-30) \put(10,-30){\vector(1,-3){0}}
 \put(10,-30){\makebox(0,0)[t]{$S_{1}$}}
 \put(-20,-20){\makebox(0,0)[rt]{$[+1]$}}
 \arc{28.284}{1.249}{2.356} \put(4.472,-13.416){\vector(3,1){0}}
 \arc{42.426}{1.249}{3.927} \put(6.708,-20.125){\vector(3,1){0}}
 \arc{14.142}{3.927}{5.498} \put(-5,5){\vector(-1,-1){0}}
\end{picture}
&
\begin{picture}(80,80)(-40,-40)
 \dottedline{3}(-20,-20)(0,0)(20,-20)%
 \put(0,0){\vector(1,1){20}} \put(20,20){\makebox(0,0)[lb]{$S_{123}$}}
 \dashline[50]{5}(0,0)(-20,20) \put(-20,20){\vector(-1,1){0}}
 \put(-20,20){\makebox(0,0)[rb]{$S_{1}$}}
 \dottedline{2}(0,0)(-10,30) \put(-10,30){\vector(-1,3){0}}
 \put(-10,30){\makebox(0,0)[b]{$S_{2}$}}
 \put(20,-20){\makebox(0,0)[lt]{$[-1]$}}
 \arc{28.284}{4.390}{5.498} \put(-4.472,13.416){\vector(-3,-1){0}}
 \arc{42.426}{4.390}{7.069} \put(-6.708,20.125){\vector(-3,-1){0}}
 \arc{14.142}{3.927}{5.498} \put(-5,5){\vector(-1,-1){0}}
\end{picture}
&
\begin{picture}(80,80)(-40,-40)
 \dottedline{3}(-20,-20)(0,0)(20,-20)%
 \put(0,0){\vector(-1,1){20}} \put(-20,20){\makebox(0,0)[rb]{$S_1$}}
 \dashline[50]{5}(0,0)(20,20) \put(20,20){\vector(1,1){0}}
 \put(20,20){\makebox(0,0)[lb]{$S_{23}$}}
 \dottedline{2}(0,0)(-30,-10) \put(-30,-10){\vector(-3,-1){0}}
 \put(-30,-10){\makebox(0,0)[r]{$S_{2}$}}
 \put(20,-20){\makebox(0,0)[lt]{$[-1]$}}
 \arc{28.284}{2.820}{3.927} \put(-13.416,-4.472){\vector(1,-3){0}}
 \arc{42.426}{2.820}{5.498} \put(-20.125,-6.708){\vector(1,-3){0}}
 \arc{14.142}{5.498}{7.069} \put(5,5){\vector(-1,1){0}}
\end{picture}
\\\hline
$I=(S_{12},S_1,S_3)$
& $J=(S_{23},S_2,S_{123})$
& $K=(S_2,S_{123},S_{12})$
& $L=(S_3,S_1,S_{23})$
\\
\begin{picture}(80,80)(-40,-40)
 \dottedline{3}(-20,-20)(0,0)(20,-20)%
 \put(0,0){\vector(1,1){20}} \put(20,20){\makebox(0,0)[lb]{$S_{12}$}}
 \dashline[50]{5}(0,0)(-20,20) \put(-20,20){\vector(-1,1){0}}
 \put(-20,20){\makebox(0,0)[rb]{$S_{1}$}}
 \dottedline{2}(0,0)(10,-30) \put(10,-30){\vector(1,-3){0}}
 \put(10,-30){\makebox(0,0)[t]{$S_{3}$}}
 \put(-20,-20){\makebox(0,0)[rt]{$[-1]$}}
 \arc{28.284}{1.249}{2.356} \put(4.472,-13.416){\vector(3,1){0}}
 \arc{14.142}{3.927}{5.498} \put(-5,5){\vector(-1,-1){0}}
\end{picture}
&
\begin{picture}(80,80)(-40,-40)
 \dottedline{3}(-20,-20)(0,0)(20,-20)%
 \put(0,0){\vector(1,1){20}} \put(20,20){\makebox(0,0)[lb]{$S_{23}$}}
 \dashline[50]{5}(0,0)(-20,20) \put(-20,20){\vector(-1,1){0}}
 \put(-20,20){\makebox(0,0)[rb]{$S_{2}$}}
 \dottedline{2}(0,0)(-10,30) \put(-10,30){\vector(-1,3){0}}
 \put(-10,30){\makebox(0,0)[b]{$S_{123}$}}
 \arc{28.284}{4.390}{5.498} \put(-4.472,13.416){\vector(-3,-1){0}}
 \arc{14.142}{3.927}{5.498} \put(-5,5){\vector(-1,-1){0}}
\end{picture}
&
\begin{picture}(80,80)(-40,-40)
 \dottedline{3}(-20,-20)(0,0)(20,-20)%
 \put(0,0){\vector(1,1){20}} \put(20,20){\makebox(0,0)[lb]{$S_{2}$}}
 \dashline[50]{5}(0,0)(-20,20) \put(-20,20){\vector(-1,1){0}}
 \put(-20,20){\makebox(0,0)[rb]{$S_{123}$}}
 \dottedline{2}(0,0)(-30,-10) \put(-30,-10){\vector(-3,-1){0}}
 \put(-30,-10){\makebox(0,0)[r]{$S_{12}$}}
 \arc{28.284}{2.820}{3.927} \put(-13.416,-4.472){\vector(1,-3){0}}
 \arc{42.426}{2.820}{5.498} \put(-20.125,-6.708){\vector(1,-3){0}}
\end{picture}
&
\begin{picture}(80,80)(-40,-40)
 \dottedline{3}(-20,-20)(0,0)(20,-20)%
 \put(0,0){\vector(1,1){20}} \put(20,20){\makebox(0,0)[lb]{$S_3$}}
 \dashline[50]{5}(0,0)(-20,20) \put(-20,20){\vector(-1,1){0}}
 \put(-20,20){\makebox(0,0)[rb]{$S_{1}$}}
 \dottedline{2}(0,0)(-10,30) \put(-10,30){\vector(-1,3){0}}
 \put(-10,30){\makebox(0,0)[b]{$S_{23}$}}
 \put(20,-20){\makebox(0,0)[lt]{$[-1]$}}
 \arc{28.284}{4.390}{5.498} \put(-4.472,13.416){\vector(-3,-1){0}}
 \arc{42.426}{4.390}{7.069} \put(-6.708,20.125){\vector(-3,-1){0}}
\end{picture}
\\\hline
\end{tabular}
\end{table}
For example, we have:
\[
 \Theta_A \cong
 \left\{\,(m_1,m_2,m_3,\phi_1,\phi_2,\phi_3)\in\R^6\,\left|\,
 \begin{array}{l}
  m_i>0,\ \phi_1<\phi_2+1,\\
  \phi_2<\phi_3+1,\ \phi_1<\phi_3+1
 \end{array}
 \,\right.\right\}
\]
\[
 m_i=|Z(S_i)|,\ \phi_i=\phi_\sigma(S_i)\ (=S_i \text{ for simply})
\]
It is easy to show that $\Theta_X=\Theta_Y$ if $X$ is equivalent
(in the sense of Araya \cite{ara}) to $Y$.

Since $D^b(A_3)$ is constructible, we denote $\Sigma(A_3)$
to be the union of all $\Theta_X$ ($X=A,\dots,L$).
It is an open and connected 3-dimensional submanifold,
and is also a connected component of $\stab(A_3)$
(we will see the proof in Appendix).

\subsection{Main Theorem}

Let $L_i=\{\,Z\in\hom_\Z(K(A_3),\C)\,|\,Z(S_i)=0\,\}$ ($i=1,2,3$),
$L_4=\{\,Z(S_{12})=0\,\}$, $L_5\{\,Z(S_{23})=0\,\}$ and
$L_6=\{\,Z(S_{123})=0\,\}$.
The subject of this paper is that the local homeomorphism
$\z:\Sigma(A_3)\rightarrow\hom_\Z(K(A_3),\C)$
becomes a covering map if it is restricted to the inverse image
of $\hom_\Z(K(A_3),\C)\setminus(L_1\cup\cdots\cup L_6)$.

\subsection{The image of the local homeomorphism}

First we show that the image of
$\z:\Sigma(A_3)\rightarrow\hom_\Z(K(A_3),\C)$ is
$\hom_\Z(K(A_3),\C)\setminus\{0\}$.
It is clear that $\mathrm{Im}(\z)$ does not contain $0$.

When a central charge $Z:K(A_3)\rightarrow\C$ vanishes at most two
$S_i$ ($i=1,2,3$), there exists a stability condition $(Z,X)$ by
appropriate choice of phases of $S_x$ and $X=A,\dots,L$.

\begin{figure}[h]
\begin{center}
\tiny
 \begin{picture}(80,80)(-40,-40)%
  \put(-40,0){\line(1,0){80}}%
  \put(0,-40){\line(0,1){80}}%
  \put(0,0){\vector(-1,-1){20}}%
  \put(-20,-20){\makebox(0,0)[rt]{$S_1$}}%
  \put(0,0){\vector(1,1){20}}%
  \put(20,20){\makebox(0,0)[lb]{$S_2$}}%
  \put(0,0){\vector(-3,1){30}}%
  \put(-40,15){\makebox(0,0)[l]{$S_3=S_{123}$}}%
  \put(5,-5){\makebox(0,0)[lt]{$(S_{12}=0)$}}%
  \put(0,0){\vector(-1,3){10}}%
  \put(-10,30){\makebox(0,0)[b]{$S_{23}$}}%
  \arc{42.426}{5.4978}{6.2831}
  \arc{14.142}{2.3562}{6.2831}
  \arc{28.284}{3.4633}{6.2831}
  \put(20,10){\makebox(0,0)[l]{$\beta$}}%
  \put(5,13){\makebox(0,0)[lb]{$\gamma$}}%
  \put(-10,-2){\makebox(0,0)[t]{$\alpha$}}%
 \end{picture}
 \hfil
 \begin{picture}(80,80)(-40,-40)%
  \put(-40,0){\line(1,0){80}}%
  \put(0,-40){\line(0,1){80}}%
  \put(0,0){\vector(1,1){20}}%
  \put(20,20){\makebox(0,0)[lb]{$S_2=S_{12}$}}%
  \put(0,0){\vector(-3,1){30}}%
  \put(-40,15){\makebox(0,0)[l]{$S_3$}}%
  \put(5,-5){\makebox(0,0)[lt]{$(S_1=0)$}}%
  \put(0,0){\vector(-1,3){10}}%
  \put(-10,30){\makebox(0,0)[b]{$S_{23}=S_{123}$}}%
  \arc{42.426}{5.4978}{6.2831}
  \arc{28.284}{3.4633}{6.2831}
  \put(20,10){\makebox(0,0)[l]{$\beta$}}%
  \put(5,13){\makebox(0,0)[lb]{$\gamma$}}%
 \end{picture}
 \hfil
 \begin{picture}(80,80)(-40,-40)%
  \put(-40,0){\line(1,0){80}}%
  \put(0,-30){\line(0,1){70}}%
  \put(0,0){\vector(-3,1){30}}%
  \put(-40,23){\makebox(0,0)[l]{$S_3=S_{23}$}}%
  \put(-29,15){\makebox(0,0)[l]{$=S_{123}$}}%
  \put(0,-40){\makebox(0,0)[b]{$(S_1=S_2=S_{12}=0)$}}%
  \arc{14.142}{3.4633}{6.2831}
  \put(5,5){\makebox(0,0)[lb]{$\gamma$}}%
 \end{picture}
\end{center}
\caption{}
\label{ex:image}
\end{figure}

For example, let $Z$ be as the left one in Figure \ref{ex:image}.
Taking phases as $S_1=\alpha$, $S_2=\beta+2$ and $S_3=\gamma+2$,
$(Z,A_{(-1,-2,-2)})$ is in $\Theta_A$.
Remark that $(Z,A_{(-1,-2,-2)})$ is a stability condition
given by the pair of a stability function and an additive full
subcategory generated by Ext-exceptional collection
$(S_1[-1],S_2[-2],S_3[-2])$.
Otherwise, taking $S_2=\beta$, $S_3=\gamma$ and $S_{123}=\gamma+2$,
$(Z,C_{(0,0,-2)})$ is in $\Theta_C$.

Let $Z$ be as the center one in Figure \ref{ex:image}.
Taking phases as $S_2=\beta$, $S_{12}=\beta+2$ and $S_3=\gamma+2$,
$(Z,B_{(0,-2,-2)})$ is in $\Theta_B$.

Let $Z$ be as the right one in Figure \ref{ex:image}.
Taking phases as $S_3=\gamma$, $S_{23}=\gamma+2$ and
$S_{123}=\gamma+4$, $(Z,D_{(0,-2,-4)})$ is in $\Theta_D$.

\subsection{Proof of the covering map property}

Next we show that if $Z$ is in $L_i$ ($i=1,\dots,6$) then any
open neighborhood $U$ of $Z$ has a non-homeomorphic inverse image.
Although, it is almost clear.
For instance, when $Z$ satisfies $Z(S_1)=0$, $Z(S_2)\neq0$ and
$Z(S_3)\neq0$, $\z^{-1}(U)\cap\Theta_A$ has a connected component
$V$ such that $\z(V)$ covers over $U$ except $Z$ many times.

Finally, we finish the proof of Main Theorem.

Suppose that $Z$ does not vanish all $S_x$ ($x\in\{1,2,3,12,23,123\}$)
and $Z(S_x)\neq\lambda Z(S_y)$ for all nonzero real number $\lambda$.
In this case, it is easy to see that $\z^{-1}(U)\cap\Theta_X$ consists
of distinct open sets, which homeomorphic to $U$, for any $X=A,\dots,L$.

The case we should concern is the neighborhood $U$ of $Z$ which
lies on a boundary of $\Theta_X$.
For example, $Z$ with $Z(S_1)=1$, $Z(S_2)=-2$ and $Z(S_3)=i$
lie on the boundary of $\Theta_A$, that is, there is a connected
component $V_A$ of $\z^{-1}(U)\cap\Theta_A$ such that it does not
surjective on $U$.
Table \ref{CovProp} show that each boundary of
$\Theta_X$ is included in the interior of some $\Theta_Y$.
For example, the inverse image of
$\{Z(S_2)=-\lambda Z(S_1)\,(\lambda>0),\,Z(S_3)\neq0\}$
is a boundary of $\Theta_A$ when we define the phase of $S_2$
is equal to the phase of $S_1$ minus $1$ and the phase of $S_3$
is grater than it, and it also is a boundary of $\Theta_B$.
However it is included in the interior of $\Theta_I$.

\begin{table}[h]
 \caption{} \label{CovProp}
\tiny
\begin{tabular}{|*{4}{c|}}
\hline
\begin{picture}(80,60)(-40,-30)
 \put(0,0){\vector(1,0){20}}%
 \put(20,-5){\makebox(0,0)[c]{$S_1$}}%
 \put(0,0){\vector(-1,0){35}}%
 \put(-35,-5){\makebox(0,0)[c]{$S_2$}}%
 \put(0,0){\vector(-1,0){15}}%
 \put(-15,-5){\makebox(0,0)[c]{$S_{12}$}}%
 \dottedline{2}(0,0)(0,-20)%
 \put(0,-20){\vector(0,-1){0}}%
 \put(0,-25){\makebox(0,0)[c]{$S_3$}}%
 \put(5,-25){\makebox(0,0)[l]{(A,B) [I]}}%
 \thicklines%
  \arc{10}{0}{3.14159}%
 \thinlines%
 \arc{20}{1.5702}{3.1415}%
 \put(0,-10){\vector(4,-1){0}}%
 \put(-20,5){\makebox(0,0)[b]{[I]}}%
\end{picture}%
&
\begin{picture}(80,60)(-40,-30)
 \put(0,0){\vector(1,0){35}}%
 \put(35,-5){\makebox(0,0)[c]{$S_1$}}%
 \put(0,0){\vector(-1,0){20}}%
 \put(-20,-5){\makebox(0,0)[c]{$S_2$}}%
 \put(0,0){\vector(1,0){15}}%
 \put(15,-5){\makebox(0,0)[c]{$S_{12}$}}%
 \dottedline{2}(0,0)(0,-20)%
 \put(0,-20){\vector(0,-1){0}}%
 \put(0,-25){\makebox(0,0)[c]{$S_3$}}%
 \put(5,-25){\makebox(0,0)[l]{(A,I) [B]}}%
 \thicklines%
  \arc{10}{0}{3.14159}%
 \thinlines%
 \arc{20}{1.5702}{3.1415}%
 \put(0,-10){\vector(4,-1){0}}%
 \dashline[10]{2}(-15,0)(-15,10)%
 \put(-15,15){\makebox(0,0)[c]{$-S_{12}$}}%
 \put(-13,8){\makebox(0,0)[l]{[C]}}%
 \put(-17,8){\makebox(0,0)[r]{[K]}}%
\end{picture}
&
\begin{picture}(80,60)(-40,-30)
 \put(0,0){\vector(1,0){20}}%
 \put(20,-5){\makebox(0,0)[c]{$S_1$}}%
 \put(0,0){\vector(1,0){35}}%
 \put(35,-5){\makebox(0,0)[c]{$S_{12}$}}%
 \put(0,0){\vector(1,0){15}}%
 \put(15,5){\makebox(0,0)[c]{$S_2$}}%
 \dottedline{2}(0,0)(0,20)%
 \put(0,20){\vector(0,1){0}}%
 \put(0,25){\makebox(0,0)[c]{$S_3$}}%
 \put(5,25){\makebox(0,0)[l]{(B,I) [A]}}%
 \dottedline{2}(0,0)(0,-20)%
 \put(0,-20){\vector(0,-1){0}}%
 \put(0,-25){\makebox(0,0)[c]{$S_3$}}%
 \put(5,-25){\makebox(0,0)[l]{(I) [A]}}%
 \dottedline{2}(-40,0)(0,0)%
 \arc{20}{4.7117}{6.2832}%
 \put(0,10){\vector(-4,1){0}}%
 \arc{20}{0}{3.1415}%
 \put(0,-10){\vector(4,-1){0}}%
 \put(10,0){\vector(1,4){0}}%
 \dashline[10]{2}(-35,0)(-35,10)%
 \dashline[10]{2}(-15,0)(-15,10)%
 \put(-35,15){\makebox(0,0)[c]{$-S_{12}$}}%
 \put(-15,15){\makebox(0,0)[c]{$-S_2$}}%
 \put(-13,8){\makebox(0,0)[l]{[L]}}%
 \put(-25,8){\makebox(0,0)[c]{[J]}}%
 \put(-35,8){\makebox(0,0)[r]{[G]}}%
\end{picture}%
&
\begin{picture}(80,60)(-40,-30)
 \put(0,0){\vector(1,0){30}}%
 \put(30,5){\makebox(0,0)[c]{$S_1$}}%
 \put(0,0){\vector(-1,0){30}}%
 \put(-30,5){\makebox(0,0)[c]{$S_3$}}%
 \dottedline{2}(0,0)(0,-20)%
 \put(0,-20){\vector(0,-1){0}}%
 \put(0,-25){\makebox(0,0)[c]{$S_2$}}%
 \put(5,-25){\makebox(0,0)[l]{(A) [I]}}%
 \thicklines
  \arc{10}{0}{3.14159}%
 \thinlines
 \arc{20}{0}{3.1415}%
 \put(10,0){\vector(1,4){0}}%
 \put(0,-10){\vector(4,-1){0}}%
\end{picture}
\\\hline
\begin{picture}(80,60)(-40,-30)
 \put(0,0){\vector(1,0){20}}%
 \put(20,5){\makebox(0,0)[c]{$S_{12}$}}%
 \put(0,0){\vector(-1,0){35}}%
 \put(-35,5){\makebox(0,0)[c]{$S_3$}}%
 \put(0,0){\vector(-1,0){15}}%
 \put(-18,5){\makebox(0,0)[c]{$S_{123}$}}%
 \dottedline{2}(0,0)(0,20)%
 \put(0,20){\vector(0,1){0}}%
 \put(0,25){\makebox(0,0)[c]{$S_2$}}%
 \put(5,25){\makebox(0,0)[l]{(B,C) [K]}}%
 \thicklines%
  \arc{10}{0}{3.14159}%
 \thinlines%
 \arc{20}{3.1415}{4.7117}%
 \put(-10,0){\vector(-1,-4){0}}%
 \put(-20,-8){\makebox(0,0)[c]{[K]}}%
\end{picture}%
&
\begin{picture}(80,60)(-40,-30)
 \put(0,0){\vector(1,0){35}}%
 \put(35,-5){\makebox(0,0)[c]{$S_{12}$}}%
 \put(0,0){\vector(-1,0){20}}%
 \put(-20,-5){\makebox(0,0)[c]{$S_3$}}%
 \put(0,0){\vector(1,0){15}}%
 \put(12,-5){\makebox(0,0)[l]{$S_{123}$}}%
 \dottedline{2}(0,0)(0,20)%
 \put(0,20){\vector(0,1){0}}%
 \put(0,25){\makebox(0,0)[c]{$S_2$}}%
 \put(5,25){\makebox(0,0)[l]{(B,K) [C]}}%
 \dottedline{2}(0,0)(0,-20)%
 \put(0,-20){\vector(0,-1){0}}%
 \put(0,-25){\makebox(0,0)[c]{$S_2$}}%
 \put(5,-25){\makebox(0,0)[l]{(K) [C]}}%
 \dottedline{2}(35,0)(40,0)%
 \thicklines%
  \arc{10}{0}{3.14159}%
 \thinlines%
 \arc{20}{3.1415}{4.7117}%
 \put(-10,0){\vector(-1,-4){0}}%
 \arc{20}{0}{3.1415}%
 \put(0,-10){\vector(4,-1){0}}%
 \put(10,0){\vector(1,4){0}}%
 \dashline[10]{2}(20,0)(20,10)%
 \dashline[10]{2}(35,0)(35,10)%
 \put(20,15){\makebox(0,0)[c]{$-S_3$}}%
 \put(19,8){\makebox(0,0)[r]{[J]}}%
 \put(28,8){\makebox(0,0)[c]{[L]}}%
 \put(36,8){\makebox(0,0)[l]{[E]}}%
\end{picture}
&
\begin{picture}(80,60)(-40,-30)
 \put(0,0){\vector(1,0){20}}%
 \put(20,-5){\makebox(0,0)[c]{$S_{12}$}}%
 \dottedline{2}(-30,0)(0,0)%
 \put(0,0){\vector(1,0){35}}%
 \put(35,-5){\makebox(0,0)[c]{$S_{123}$}}%
 \put(0,0){\vector(1,0){15}}%
 \put(15,-10){\makebox(0,0)[c]{$S_3$}}%
 \dottedline{2}(0,0)(0,-20)%
 \put(0,-20){\vector(0,-1){0}}%
 \put(0,-25){\makebox(0,0)[c]{$S_2$}}%
 \put(5,-25){\makebox(0,0)[l]{(C,K) [B]}}%
 \arc{20}{0}{1.5707}%
 \put(10,0){\vector(1,4){0}}%
 \dashline[10]{2}(20,0)(20,10)%
 \put(19,8){\makebox(0,0)[r]{[A]}}%
 \put(21,8){\makebox(0,0)[l]{[I]}}%
\end{picture}%
&
\begin{picture}(80,60)(-40,-30)
 \put(0,0){\vector(1,0){30}}%
 \put(30,5){\makebox(0,0)[c]{$S_2$}}%
 \dottedline{2}(0,0)(-30,0)%
 \put(0,0){\vector(1,0){20}}%
 \put(20,-5){\makebox(0,0)[c]{$S_3$}}%
 \dottedline{2}(0,0)(0,20)%
 \put(0,20){\vector(0,1){0}}%
 \put(0,25){\makebox(0,0)[c]{$S_{12}$}}%
 \put(8,25){\makebox(0,0)[l]{(B) [K]}}%
 \arc{20}{3.1415}{6.2832}%
 \put(-10,0){\vector(-1,-4){0}}%
 \put(0,10){\vector(-4,1){0}}%
\end{picture}
\\\hline
\begin{picture}(80,60)(-40,-30)
 \put(0,0){\vector(1,0){20}}%
 \put(20,5){\makebox(0,0)[c]{$S_2$}}%
 \put(0,0){\vector(-1,0){35}}%
 \put(-35,5){\makebox(0,0)[c]{$S_3$}}%
 \put(0,0){\vector(-1,0){15}}%
 \put(-15,5){\makebox(0,0)[c]{$S_{23}$}}%
 \dottedline{2}(0,0)(0,20)%
 \put(0,20){\vector(0,1){0}}%
 \put(0,25){\makebox(0,0)[c]{$S_{123}$}}%
 \put(10,25){\makebox(0,0)[l]{(C,D) [J]}}%
 \thicklines%
  \arc{10}{0}{3.14159}%
 \thinlines%
 \arc{20}{4.7117}{6.2832}%
 \put(0,10){\vector(-4,1){0}}%
 \put(20,-8){\makebox(0,0)[c]{[J]}}%
\end{picture}%
&
\begin{picture}(80,60)(-40,-30)
 \put(0,0){\vector(1,0){35}}%
 \put(35,5){\makebox(0,0)[c]{$S_2$}}%
 \put(0,0){\vector(-1,0){15}}%
 \put(-15,5){\makebox(0,0)[c]{$S_3$}}%
 \put(0,0){\vector(1,0){20}}%
 \put(20,5){\makebox(0,0)[c]{$S_{23}$}}%
 \dottedline{2}(0,0)(0,20)%
 \put(0,20){\vector(0,1){0}}%
 \put(0,25){\makebox(0,0)[c]{$S_{123}$}}%
 \put(10,25){\makebox(0,0)[l]{(C,J) [D]}}%
 \thicklines
  \arc{10}{0}{3.14159}%
 \thinlines
 \arc{20}{4.7117}{6.2832}%
 \put(0,10){\vector(-4,1){0}}%
 \dashline[10]{2}(20,0)(20,-10)%
 \put(19,-8){\makebox(0,0)[r]{[E]}}%
 \put(21,-8){\makebox(0,0)[l]{[L]}}%
\end{picture}
&
\begin{picture}(80,60)(-40,-30)
 \put(0,0){\vector(1,0){20}}%
 \put(20,5){\makebox(0,0)[l]{$S_2$}}%
 \dottedline{2}(-30,0)(0,0)%
 \put(0,0){\vector(1,0){35}}%
 \put(35,5){\makebox(0,0)[c]{$S_{23}$}}%
 \put(0,0){\vector(1,0){15}}%
 \put(15,5){\makebox(0,0)[c]{$S_3$}}%
 \dottedline{2}(0,0)(0,-20)%
 \put(0,-20){\vector(0,-1){0}}%
 \put(0,-25){\makebox(0,0)[c]{$S_{123}$}}%
 \put(10,25){\makebox(0,0)[l]{(J) [C]}}%
 \dottedline{2}(0,0)(0,20)%
 \put(0,20){\vector(0,1){0}}%
 \put(0,25){\makebox(0,0)[c]{$S_{123}$}}%
 \put(10,-25){\makebox(0,0)[l]{(D,J) [C]}}%
 \dottedline{2}(35,0)(40,0)%
 \arc{20}{3.1415}{6.2832}%
 \put(0,10){\vector(-4,1){0}}%
 \put(-10,0){\vector(-1,-4){0}}%
 \arc{20}{1.5707}{3.1415}%
 \put(0,-10){\vector(4,-1){0}}%
 \dashline[10]{2}(15,0)(15,-10)%
 \dashline[10]{2}(35,0)(35,-10)%
 \put(15,-8){\makebox(0,0)[r]{[K]}}%
 \put(25,-8){\makebox(0,0)[c]{[I]}}%
 \put(35,-8){\makebox(0,0)[l]{[A]}}%
\end{picture}%
&
\begin{picture}(80,60)(-40,-30)
 \put(0,0){\vector(1,0){30}}%
 \put(30,5){\makebox(0,0)[c]{$S_2$}}%
 \dottedline{2}(0,0)(-30,0)%
 \put(0,0){\vector(1,0){20}}%
 \put(20,-5){\makebox(0,0)[c]{$S_{123}$}}%
 \dottedline{2}(0,0)(0,-20)%
 \put(0,-20){\vector(0,-1){0}}%
 \put(0,-25){\makebox(0,0)[c]{$S_3$}}%
 \put(5,-25){\makebox(0,0)[l]{(C) [J]}}%
 \arc{20}{0}{3.1415}%
 \put(10,0){\vector(1,4){0}}%
 \put(0,-10){\vector(4,-1){0}}%
\end{picture}
\\\hline
\begin{picture}(80,60)(-40,-30)
 \put(0,0){\vector(1,0){20}}%
 \put(20,-5){\makebox(0,0)[c]{$S_{23}$}}%
 \dottedline{2}(-30,0)(0,0)%
 \put(0,0){\vector(1,0){35}}%
 \put(35,-5){\makebox(0,0)[c]{$S_{123}$}}%
 \put(0,0){\vector(1,0){15}}%
 \put(15,5){\makebox(0,0)[c]{$S_1$}}%
 \dottedline{2}(0,0)(0,20)%
 \put(0,20){\vector(0,1){0}}%
 \put(0,25){\makebox(0,0)[c]{$S_3$}}%
 \put(5,25){\makebox(0,0)[l]{(D,E) [L]}}%
 \thicklines%
  \arc{10}{0}{3.14159}%
 \thinlines%
 \arc{20}{3.1415}{4.7117}%
 \put(-10,0){\vector(-1,-4){0}}%
 \put(-20,-8){\makebox(0,0)[c]{[L]}}%
\end{picture}%
&
\begin{picture}(80,60)(-40,-30)
 \put(0,0){\vector(1,0){35}}%
 \put(35,-5){\makebox(0,0)[c]{$S_{23}$}}%
 \put(0,0){\vector(1,0){20}}%
 \put(20,-5){\makebox(0,0)[c]{$S_{123}$}}%
 \put(0,0){\vector(-1,0){15}}%
 \put(-15,-5){\makebox(0,0)[c]{$S_1$}}%
 \dottedline{2}(0,0)(0,20)%
 \put(0,20){\vector(0,1){0}}%
 \put(0,25){\makebox(0,0)[c]{$S_3$}}%
 \put(5,25){\makebox(0,0)[l]{(D,L) [E]}}%
 \dottedline{2}(0,0)(0,-20)%
 \put(0,-20){\vector(0,-1){0}}%
 \put(0,-25){\makebox(0,0)[c]{$S_3$}}%
 \put(5,-25){\makebox(0,0)[l]{(L) [E]}}%
 \dottedline{2}(35,0)(40,0)%
 \thicklines%
  \arc{10}{3.14159}{6.2832}%
 \thinlines%
 \arc{20}{0}{4.7117}%
 \put(10,0){\vector(1,4){0}}%
 \put(-10,0){\vector(-1,-4){0}}%
 \put(0,-10){\vector(4,-1){0}}%
 \dashline[10]{2}(20,0)(20,10)%
 \dashline[10]{2}(35,0)(35,10)%
 \put(19,8){\makebox(0,0)[r]{[I]}}%
 \put(28,8){\makebox(0,0)[c]{[K]}}%
 \put(35,8){\makebox(0,0)[l]{[G]}}%
\end{picture}
&
\begin{picture}(80,60)(-40,-30)
 \put(0,0){\vector(1,0){20}}%
 \put(20,-5){\makebox(0,0)[c]{$S_{23}$}}%
 \put(0,0){\vector(-1,0){35}}%
 \put(-35,-5){\makebox(0,0)[c]{$S_1$}}%
 \put(0,0){\vector(-1,0){15}}%
 \put(-15,-5){\makebox(0,0)[c]{$S_{123}$}}%
 \dottedline{2}(0,0)(0,-20)%
 \put(0,-20){\vector(0,-1){0}}%
 \put(0,-25){\makebox(0,0)[c]{$S_3$}}%
 \put(5,-25){\makebox(0,0)[l]{(E,L) [D]}}%
 \dottedline{2}(20,0)(35,0)%
 \thicklines%
  \arc{10}{3.1415}{6.2832}%
 \thinlines%
 \arc{20}{0}{1.5707}%
 \put(10,0){\vector(1,4){0}}%
 \dashline[10]{2}(20,0)(20,10)%
 \put(19,8){\makebox(0,0)[r]{[C]}}%
 \put(21,8){\makebox(0,0)[l]{[J]}}%
\end{picture}%
&
\begin{picture}(80,60)(-40,-30)
 \put(0,0){\vector(1,0){30}}%
 \put(30,5){\makebox(0,0)[c]{$S_{123}$}}%
 \put(0,0){\vector(-1,0){20}}%
 \put(-20,-5){\makebox(0,0)[c]{$S_3$}}%
 \dottedline{2}(0,0)(0,-20)%
 \put(0,-20){\vector(0,-1){0}}%
 \put(0,-25){\makebox(0,0)[c]{$S_{23}$}}%
 \put(8,-25){\makebox(0,0)[l]{(D) [L]}}%
 \thicklines%
  \arc{10}{0}{3.14159}%
 \thinlines%
 \arc{20}{0}{3.1415}%
 \put(10,0){\vector(1,4){0}}%
 \put(0,-10){\vector(4,-1){0}}%
\end{picture}
\\\hline
\begin{picture}(80,60)(-40,-30)
 \put(0,0){\vector(-1,0){20}}%
 \put(-20,5){\makebox(0,0)[c]{$S_3$}}%
 \put(0,0){\vector(-1,0){35}}%
 \put(-35,5){\makebox(0,0)[c]{$S_{123}$}}%
 \put(0,0){\vector(-1,0){15}}%
 \put(-15,5){\makebox(0,0)[l]{$S_{12}$}}%
 \dottedline{2}(0,0)(0,20)%
 \put(0,20){\vector(0,1){0}}%
 \put(0,25){\makebox(0,0)[c]{$S_1$}}%
 \put(5,25){\makebox(0,0)[l]{(E,F) [I]}}%
 \dottedline{2}(0,0)(35,0)%
 \thicklines%
  \arc{10}{0}{3.1415}%
 \thinlines%
 \arc{20}{4.7123}{6.2832}%
 \put(0,10){\vector(-4,1){0}}%
 \put(20,-8){\makebox(0,0)[c]{[I]}}%
\end{picture}%
&
\begin{picture}(80,60)(-40,-30)
 \put(0,0){\vector(-1,0){35}}%
 \put(-35,5){\makebox(0,0)[c]{$S_3$}}%
 \put(0,0){\vector(1,0){20}}%
 \put(20,5){\makebox(0,0)[c]{$S_{12}$}}%
 \put(0,0){\vector(-1,0){15}}%
 \put(-15,5){\makebox(0,0)[c]{$S_{123}$}}%
 \dottedline{2}(0,0)(0,20)%
 \put(0,20){\vector(0,1){0}}%
 \put(0,25){\makebox(0,0)[c]{$S_1$}}%
 \put(5,25){\makebox(0,0)[l]{(E,I) [F]}}%
 \dottedline{2}(20,0)(35,0)%
 \thicklines%
  \arc{10}{0}{3.14159}%
 \thinlines%
 \arc{20}{4.7123}{6.2832}%
 \put(0,10){\vector(-4,1){0}}%
 \dashline[10]{2}(20,0)(20,-10)%
 \put(18,-8){\makebox(0,0)[r]{[G]}}%
 \put(22,-8){\makebox(0,0)[l]{[K]}}%
\end{picture}
&
\begin{picture}(80,60)(-40,-30)
 \put(0,0){\vector(-1,0){15}}%
 \put(-15,5){\makebox(0,0)[r]{$S_3$}}%
 \put(0,0){\vector(1,0){35}}%
 \put(35,5){\makebox(0,0)[c]{$S_{12}$}}%
 \put(0,0){\vector(1,0){20}}%
 \put(12,5){\makebox(0,0)[l]{$S_{123}$}}%
 \dottedline{2}(0,0)(0,-20)%
 \put(0,-20){\vector(0,-1){0}}%
 \put(0,-25){\makebox(0,0)[c]{$S_1$}}%
 \put(5,-25){\makebox(0,0)[l]{(F,I) [E]}}%
 \dottedline{2}(0,0)(0,20)%
 \put(0,20){\vector(0,1){0}}%
 \put(0,25){\makebox(0,0)[c]{$S_1$}}%
 \put(5,25){\makebox(0,0)[l]{(I) [E]}}%
 \dottedline{2}(35,0)(40,0)%
 \thicklines%
  \arc{10}{0}{3.14159}%
 \thinlines%
 \arc{20}{1.5707}{3.1415}%
 \put(0,-10){\vector(4,-1){0}}%
 \arc{20}{3.1415}{6.2832}%
 \put(0,10){\vector(-4,1){0}}%
 \put(-10,0){\vector(-1,-4){0}}%
 \dashline[10]{2}(20,0)(20,-10)%
 \dashline[10]{2}(35,0)(35,-10)%
 \put(18,-8){\makebox(0,0)[r]{[J]}}%
 \put(28,-8){\makebox(0,0)[c]{[L]}}%
 \put(35,-8){\makebox(0,0)[l]{[C]}}%
\end{picture}%
&
\begin{picture}(80,60)(-40,-30)
 \put(0,0){\vector(1,0){30}}%
 \put(30,5){\makebox(0,0)[c]{$S_1$}}%
 \put(0,0){\vector(-1,0){30}}%
 \put(-30,5){\makebox(0,0)[c]{$S_3$}}%
 \dottedline{2}(0,0)(0,-20)%
 \put(0,-20){\vector(0,-1){0}}%
 \put(0,-25){\makebox(0,0)[c]{$S_{123}$}}%
 \put(10,-25){\makebox(0,0)[l]{(E) [I]}}%
 \thicklines%
  \arc{10}{0}{3.14159}%
 \thinlines%
 \arc{20}{0}{3.1415}%
 \put(10,0){\vector(1,4){0}}%
 \put(0,-10){\vector(4,-1){0}}%
\end{picture}
\\\hline
\begin{picture}(80,60)(-40,-30)
 \put(0,0){\vector(1,0){35}}%
 \put(35,5){\makebox(0,0)[c]{$S_1$}}%
 \put(0,0){\vector(1,0){15}}%
 \put(15,5){\makebox(0,0)[c]{$S_{12}$}}%
 \put(0,0){\vector(-1,0){20}}%
 \put(-20,5){\makebox(0,0)[c]{$S_2$}}%
 \dottedline{2}(0,0)(0,20)%
 \put(0,20){\vector(0,1){0}}%
 \put(0,25){\makebox(0,0)[c]{$S_{123}$}}%
 \put(10,25){\makebox(0,0)[l]{(F,G) [K]}}%
 \thicklines%
  \arc{10}{0}{3.1415}%
 \thinlines%
 \arc{20}{3.1415}{4.7123}%
 \put(-10,0){\vector(-1,-4){0}}%
 \put(-20,-8){\makebox(0,0)[c]{[K]}}%
\end{picture}
&
\begin{picture}(80,60)(-40,-30)
 \put(0,0){\vector(1,0){35}}%
 \put(35,-5){\makebox(0,0)[c]{$S_{12}$}}%
 \put(0,0){\vector(1,0){20}}%
 \put(20,-5){\makebox(0,0)[c]{$S_1$}}%
 \put(0,0){\vector(1,0){15}}%
 \put(15,-8){\makebox(0,0)[c]{$S_2$}}%
 \dottedline{2}(0,0)(0,-20)%
 \put(0,-20){\vector(0,-1){0}}%
 \put(0,-25){\makebox(0,0)[c]{$S_{123}$}}%
 \put(10,-25){\makebox(0,0)[l]{(K) [G]}}%
 \dottedline{2}(0,0)(0,20)%
 \put(0,20){\vector(0,1){0}}%
 \put(0,25){\makebox(0,0)[c]{$S_{123}$}}%
 \put(10,25){\makebox(0,0)[l]{(F,K) [G]}}%
 \dottedline{2}(0,0)(-30,0)%
 \dottedline{2}(35,0)(40,0)%
 \thicklines%
  \arc{10}{0}{3.1415}%
 \thinlines%
 \arc{20}{0}{4.7122}%
 \put(-10,0){\vector(-1,-4){0}}%
 \put(0,-10){\vector(4,-1){0}}%
 \put(10,0){\vector(1,4){0}}%
 \dashline[10]{2}(20,0)(20,10)%
 \dashline[10]{2}(35,0)(35,10)%
 \put(19,8){\makebox(0,0)[r]{[J]}}%
 \put(28,8){\makebox(0,0)[c]{[L]}}%
 \put(35,8){\makebox(0,0)[l]{[A]}}%
\end{picture}
&
\begin{picture}(80,60)(-40,-30)
 \put(0,0){\vector(1,0){15}}%
 \put(15,5){\makebox(0,0)[l]{$S_1$}}%
 \put(0,0){\vector(-1,0){35}}%
 \put(-35,5){\makebox(0,0)[c]{$S_2$}}%
 \put(0,0){\vector(-1,0){20}}%
 \put(-20,5){\makebox(0,0)[c]{$S_{12}$}}%
 \dottedline{2}(0,0)(0,20)%
 \put(0,20){\vector(0,1){0}}%
 \put(0,25){\makebox(0,0)[c]{$S_{123}$}}%
 \put(10,25){\makebox(0,0)[l]{(G,K) [F]}}%
 \thicklines%
  \arc{10}{0}{3.1415}%
 \thinlines%
 \arc{20}{3.1415}{4.7122}%
 \put(-10,0){\vector(-1,-4){0}}%
 \dashline[10]{2}(-20,0)(-20,-10)%
 \put(-19,-8){\makebox(0,0)[l]{[E]}}%
 \put(-21,-8){\makebox(0,0)[r]{[I]}}%
\end{picture}
&
\begin{picture}(80,60)(-40,-30)
 \put(0,0){\vector(-1,0){20}}%
 \put(-20,5){\makebox(0,0)[c]{$S_{123}$}}%
 \put(0,0){\vector(1,0){30}}%
 \put(30,5){\makebox(0,0)[c]{$S_1$}}%
 \dottedline{2}(0,0)(0,-20)%
 \put(0,-20){\vector(0,-1){0}}%
 \put(0,-25){\makebox(0,0)[c]{$S_{12}$}}%
 \put(8,-25){\makebox(0,0)[l]{(F) [K]}}%
 \thicklines%
  \arc{10}{0}{3.1415}%
 \thinlines%
 \arc{20}{0}{3.1415}%
 \put(0,-10){\vector(4,-1){0}}%
 \put(10,0){\vector(1,4){0}}%
\end{picture}
\\\hline
\begin{picture}(80,60)(-40,-30)
 \put(0,0){\vector(1,0){15}}%
 \put(10,5){\makebox(0,0)[l]{$S_{123}$}}%
 \put(0,0){\vector(1,0){35}}%
 \put(35,5){\makebox(0,0)[c]{$S_1$}}%
 \put(0,0){\vector(-1,0){20}}%
 \put(-20,5){\makebox(0,0)[c]{$S_{23}$}}%
 \dottedline{2}(0,0)(0,20)%
 \put(0,20){\vector(0,1){0}}%
 \put(0,25){\makebox(0,0)[c]{$S_2$}}%
 \put(5,25){\makebox(0,0)[l]{(G,H) [J]}}%
 \thicklines%
  \arc{10}{0}{3.1415}%
 \thinlines%
 \arc{20}{4.7123}{6.2831}%
 \put(0,10){\vector(-4,1){0}}%
 \put(20,-8){\makebox(0,0)[c]{[J]}}%
\end{picture}
&
\begin{picture}(80,60)(-40,-30)
 \put(0,0){\vector(1,0){35}}%
 \put(35,5){\makebox(0,0)[c]{$S_{123}$}}%
 \put(0,0){\vector(1,0){15}}%
 \put(15,10){\makebox(0,0)[r]{$S_1$}}%
 \put(0,0){\vector(1,0){20}}%
 \put(20,5){\makebox(0,0)[c]{$S_{23}$}}%
 \dottedline{2}(0,0)(0,20)%
 \put(0,20){\vector(0,1){0}}%
 \put(0,25){\makebox(0,0)[c]{$S_2$}}%
 \put(5,25){\makebox(0,0)[l]{(G,J) [H]}}%
 \dottedline{2}(0,0)(-30,0)%
 \thicklines%
  \arc{10}{0}{3.1415}%
 \thinlines%
 \arc{20}{4.7123}{6.2831}%
 \put(0,10){\vector(-4,1){0}}%
 \dashline[10]{2}(20,0)(20,-10)%
 \put(19,-8){\makebox(0,0)[r]{[A]}}%
 \put(21,-8){\makebox(0,0)[l]{[L]}}%
\end{picture}
&
\begin{picture}(80,60)(-40,-30)
 \put(0,0){\vector(-1,0){15}}%
 \put(-10,-5){\makebox(0,0)[r]{$S_{123}$}}%
 \put(0,0){\vector(1,0){20}}%
 \put(20,-5){\makebox(0,0)[l]{$S_1$}}%
 \put(0,0){\vector(-1,0){35}}%
 \put(-35,-5){\makebox(0,0)[c]{$S_{23}$}}%
 \dottedline{2}(0,0)(0,20)%
 \put(0,20){\vector(0,1){0}}%
 \put(0,25){\makebox(0,0)[c]{$S_2$}}%
 \put(5,25){\makebox(0,0)[l]{(H,J) [G]}}%
 \dottedline{2}(0,0)(0,-20)%
 \put(0,-20){\vector(0,-1){0}}%
 \put(0,-25){\makebox(0,0)[c]{$S_2$}}%
 \put(5,-25){\makebox(0,0)[l]{(J) [G]}}%
 \dottedline{2}(-40,0)(-35,0)%
 \thicklines%
  \arc{10}{0}{3.1415}%
 \thinlines%
 \arc{20}{0}{3.1415}%
 \put(0,-10){\vector(4,-1){0}}%
 \put(10,0){\vector(1,4){0}}%
 \arc{20}{4.7123}{6.2831}%
 \put(0,10){\vector(-4,1){0}}%
 \dashline[10]{2}(-35,0)(-35,10)%
 \dashline[10]{2}(-20,0)(-20,10)%
 \put(-20,15){\makebox(0,0)[c]{$-S_1$}}%
 \put(-19,8){\makebox(0,0)[l]{[K]}}%
 \put(-28,8){\makebox(0,0)[c]{[I]}}%
 \put(-35,8){\makebox(0,0)[r]{[E]}}%
\end{picture}
&
\begin{picture}(80,60)(-40,-30)
 \put(0,0){\vector(1,0){35}}%
 \put(35,5){\makebox(0,0)[c]{$S_{123}$}}%
 \put(0,0){\vector(1,0){20}}%
 \put(20,5){\makebox(0,0)[c]{$S_2$}}%
 \dottedline{2}(0,0)(0,20)%
 \put(0,20){\vector(0,1){0}}%
 \put(0,25){\makebox(0,0)[c]{$S_1$}}%
 \put(5,25){\makebox(0,0)[l]{(G) [J]}}%
 \dottedline{2}(0,0)(-30,0)%
 \arc{20}{3.1415}{6.2831}%
 \put(0,10){\vector(-4,1){0}}%
 \put(-10,0){\vector(-1,-4){0}}%
\end{picture}
\\ \hline
\begin{picture}(80,60)(-40,-30)
 \put(0,0){\vector(1,0){20}}%
 \put(20,-5){\makebox(0,0)[c]{$S_{23}$}}%
 \put(0,0){\vector(1,0){35}}%
 \put(35,-5){\makebox(0,0)[c]{$S_2$}}%
 \put(0,0){\vector(-1,0){15}}%
 \put(-15,-5){\makebox(0,0)[r]{$S_3$}}%
 \dottedline{2}(0,0)(0,-20)%
 \put(0,-20){\vector(0,-1){0}}%
 \put(0,-25){\makebox(0,0)[c]{$S_1$}}%
 \put(5,-25){\makebox(0,0)[l]{(H,A) [L]}}%
 \thicklines%
  \arc{10}{0}{3.1415}%
 \thinlines%
 \arc{20}{0}{1.5708}%
 \put(10,0){\vector(1,4){0}}%
 \put(20,8){\makebox(0,0)[c]{[L]}}%
\end{picture}
&
\begin{picture}(80,60)(-40,-30)
 \put(0,0){\vector(1,0){35}}%
 \put(35,5){\makebox(0,0)[c]{$S_{23}$}}%
 \put(0,0){\vector(1,0){20}}%
 \put(20,5){\makebox(0,0)[c]{$S_2$}}%
 \put(0,0){\vector(1,0){15}}%
 \put(15,-5){\makebox(0,0)[c]{$S_3$}}%
 \dottedline{2}(0,0)(0,20)%
 \put(0,20){\vector(0,1){0}}%
 \put(0,25){\makebox(0,0)[c]{$S_1$}}%
 \put(5,25){\makebox(0,0)[l]{(L) [A]}}%
 \dottedline{2}(0,0)(0,-20)%
 \put(0,-20){\vector(0,-1){0}}%
 \put(0,-25){\makebox(0,0)[c]{$S_1$}}%
 \put(5,-25){\makebox(0,0)[l]{(H,L) [A]}}%
 \dottedline{2}(0,0)(-35,0)%
 \thicklines%
  \arc{10}{0}{3.1415}%
 \thinlines%
 \arc{20}{0}{1.5708}%
 \put(10,0){\vector(1,4){0}}%
 \arc{20}{3.1415}{6.2832}%
 \put(0,10){\vector(-4,1){0}}%
 \put(-10,0){\vector(-1,-4){0}}%
 \dashline[10]{2}(-35,0)(-35,-10)%
 \dashline[10]{2}(-20,0)(-20,-10)%
 \put(-35,-15){\makebox(0,0)[c]{$-S_{23}$}}%
 \put(-25,-15){\makebox(0,0)[l]{$-S_2$}}%
 \put(-35,-8){\makebox(0,0)[r]{[C]}}%
 \put(-28,-8){\makebox(0,0)[c]{[K]}}%
 \put(-19,-8){\makebox(0,0)[l]{[I]}}%
\end{picture}
&
\begin{picture}(80,60)(-40,-30)
 \put(0,0){\vector(1,0){15}}%
 \put(15,-5){\makebox(0,0)[l]{$S_2$}}%
 \put(0,0){\vector(-1,0){35}}%
 \put(-35,-5){\makebox(0,0)[c]{$S_3$}}%
 \put(0,0){\vector(-1,0){20}}%
 \put(-20,-5){\makebox(0,0)[c]{$S_{23}$}}%
 \dottedline{2}(0,0)(0,-20)%
 \put(0,-20){\vector(0,-1){0}}%
 \put(0,-25){\makebox(0,0)[c]{$S_1$}}%
 \put(5,-25){\makebox(0,0)[l]{(A,L) [H]}}%
 \dottedline{2}(15,0)(35,0)%
 \thicklines%
  \arc{10}{0}{3.1415}%
 \thinlines%
 \arc{20}{0}{1.5708}%
 \put(10,0){\vector(1,4){0}}%
 \dashline[10]{2}(20,0)(20,10)%
 \put(20,15){\makebox(0,0)[c]{$-S_{23}$}}%
 \put(19,8){\makebox(0,0)[r]{[G]}}%
 \put(21,8){\makebox(0,0)[l]{[J]}}%
\end{picture}
&
\tiny
\begin{picture}(80,60)(-40,-30)
 \put(0,0){\vector(1,0){20}}%
 \put(20,5){\makebox(0,0)[c]{$S_2$}}%
 \put(0,0){\vector(1,0){35}}%
 \put(35,5){\makebox(0,0)[c]{$S_1$}}%
 \dottedline{2}(0,0)(0,-20)%
 \put(0,-20){\vector(0,-1){0}}%
 \put(0,-25){\makebox(0,0)[c]{$S_{23}$}}%
 \put(8,-25){\makebox(0,0)[l]{(H) [L]}}%
 \dottedline{2}(0,0)(-30,0)%
 \thicklines%
  \arc{10}{0}{3.1415}%
 \thinlines%
 \arc{20}{0}{3.1415}%
 \put(10,0){\vector(1,4){0}}%
 \put(0,-10){\vector(4,-1){0}}%
\end{picture}
\\\hline
\end{tabular}
\end{table}

Now we conclude that the $\z^{-1}(U)$ is a disjoint union
of connected components homeomorphic to $U$ for enough small open
neighborhood of $Z$ in $\hom_\Z(K(A_3),\C)\setminus\{L_1,\dots,L_6\}$.

\appendix
\section{}\label{sec4}

We prove that $\Sigma(A_3)$ is a connected component of $\stab(A_3)$.
It is enough to show that $\overline{\Sigma(A_3)}=\Sigma(A_3)$
since $\Sigma(A_3)$ is an open submanifold of $\stab(A_3)$.

Let $\sigma=(Z,\P)$ lie in the boundary of $\overline{\Sigma(A_3)}$.
There is a sequence $\{\sigma_i=(Z_i,\P_i)\}_{i=1,2,3,\dots}$
consisting of stability conditions in $\Sigma(A_3)$.
We can assume all $\sigma_i$ belongs to $\Theta_\E$ for some
exceptional collection $\E$ since $\Sigma(A_3)$ is a union of
finitely many $\Theta_\E$'s.
Put $\E=(E_1,E_2,E_3)$.
Then all $E_j$'s are $\sigma_i$-stable for all $i$.

When we denote $X$ to be a subset of $\stab(A_3)$ consisting of all
stability conditions in which $E_i$ ($i=1,2,3$) are semistable,
$X$ is a closed subset because of a subset consisting of all
stability conditions which has an object as semistable is closed 
(\cite[Section~5]{bri}).
Thus $\Theta_\E$ is included in $X$ and
$\overline{\Theta_\E}\subset \overline{X}=X$.
Therefore $\sigma$ is included in $X$, so $\sigma$ has all $E_i$'s
as semistable.

Now we have $Z(E_i)\neq 0$ ($i=1,2,3$) and phases $\phi_\sigma(E_i)$
are defined.
Because $\sigma$ lies in the boundary of $\Theta_\E$,
$\phi_\sigma(E_i)=\phi_\sigma(E_j)+\alpha_{i,j}$ for one or two
of $(i,j)=(1,2)$, $(1,3)$, or $(2,3)$ and
$m_\sigma(E_i)\neq m_\sigma(E_j)$ if $\alpha_{i,j}=k_{i,j}$.

However all possible stability conditions which belongs to the
boundary of $\Theta_\E$ are listed in Table~\ref{CovProp} and
it shows each such stability condition is included in some
$\Theta_\mathcal{F}$.
This shows that $\overline{\Sigma(A_3)}\subset\Sigma(A_3)$.


\end{document}